\documentclass{amsart}

\usepackage{amssymb}

\ifx\pdfoutput\undefined

\usepackage{colordvi}

\else

\input pdfcolor.tex

\fi

   \def\okr{\stackrel{\scriptscriptstyle{\sf{def}}}{=}}

\newtheorem{thm}{Theorem}[section]
\newtheorem{cor}[thm]{Corollary}

\newtheorem{prop}[thm]{Proposition}

\theoremstyle{definition}

\newtheorem{rem}[thm]{Remark}
\newtheorem{exa}[thm]{Example}

\newcommand{\Bibitem}[1]{\bibitem{#1}}

\numberwithin{equation}{section}

\def\llf{\mathfrak L}

\def\okr{\stackrel{\scriptscriptstyle{\sf{def}}}{=}}
\def\dz#1{\mathcal D({#1})}
\def\gw{^*}
\def\is#1#2{\langle#1,#2\rangle}
\def\liczp#1{{${#1}^{\text {\rm o}}$}}
\def\zb#1#2{\{{#1}\colon\ {#2}\}}

\begin{document}
\title[Adjoints and Formal Adjoints]
{Adjoints and Formal Adjoints of Matrices of Unbounded Operators}
\author{Manfred M\"oller}
\author{Franciszek Hugon Szafraniec}

\address{M. M\"oller, The John Knopfmacher Centre for Applicable
Analysis and Number Theory, School of Mathematics, University of
the Witwatersrand, WITS, 2050, South Africa, e-mail: {\tt
manfred@maths.wits.ac.za}}
\address{F. H. Szafraniec, Instytut Mathematyki, Uniwersytet
Jagiello\'nski, ul. Reymonta 4, Pl-30059 Krak\'ow, e-mail: {\tt
fhszafra@im.uj.edu.pl}}


\thanks{This work was supported by a grant of the NRF of South Africa, GUN
2053746, and by the grant KBN 2 P 03A 637 024 (Poland).}

   \begin{abstract}
   In this paper we {\em discuss} diverse aspects of mutual
   relationship between adjoints and formal adjoints of unbounded
   operators bearing a matrix structure. We emphasize on the
   behaviour of row and column operators as they turn out to be
   the germs of an arbitrary matrix operator, providing most of
   the information about the latter {as it is the troublemaker}.
   \end{abstract}

\maketitle

 \section{Introduction}
In recent years, $2\times 2$ matrices  of unbounded operators have attracted considerable
attention, roughly divided into two groups of problems{: they} occur as generators for
semigroups, see \cite{BH}, \cite{N} and \cite{Eng}, and as tools
in problems from Mathematical Physics, see \cite{HL}, \cite{F},
\cite{L}. The latter has attracted much interest in spectral
properties of such matrices, in particular in its essential
spectrum,  see 
\cite{ALMS}, 
\cite{HMM1},
\cite{KM}, \cite{KN}
and the references therein.

 Adjoints of operators bearing a matrix structure have been
investigated case by case; it might be difficult to find any
general approach to the problem when there are rows or columns
with more than one unbounded entry. So as to mention some partial
results let us refer
 to \cite{HKM} and \cite{OS} where {the only nonzero entries are those off the diagonal}.
 On the other hand, in    \cite{DGS}
  and \cite{N} examples are given showing that
{the adjoint of a matrix operator $A$ and its formal adjoint
$A^\times $, that is the matrix of adjoints of all the particular
entries,} may be quite different. This supports the idea of the
present paper to build a common framework for all the cases.
Positive results in this matter are intertwined with
counterexamples; the latters indicate that $A^\times $ may not
contain enough information on $A'$ itself. Thus, what is
essentially trivial for bounded operators appears to become
erratic for unbounded operators.

It turns out that the study of row and column operators {separately} is pretty
much helpful as a matrix can be decomposed in a sense by means of
these two. Indeed, column and row operators have a more
predictable behaviour, which helps to understand the
{2 by 2 matrix case -- this case is interesting as well as difficult enough to deserve a careful treatment}.  The pith of the problem can be simply described by saying that
   \begin{quote}
if $A$ is either a row or a column operator (even in a Hilbert space)
 with both entries being unbounded and closed then it is rather
unlikely $A$ to coincide with it second adjoint.
   \end{quote}
   On the other hand, it is important to stress that if a column
has only one entry which is unbounded then the problem does not
appear at all; {this is the most frequent case which occurs in the
literature.}

   \section{Preliminaries}
Reasonable assumptions on $2\times 2$ operator matrices
 are that their entries are closed densely defined operators and
that the operators they determine are densely defined. However,
when one divides bigger operator matrices into $2\times 2$ blocks,
the blocks in the latter structure will not be closed, in general.
Matrices of arbitrary size occur in an attempt to determine normal
extensions of unbounded operators, see \cite{S}, and our intention
is to elaborate somewhere else on this kind of matrices from the
point of view of the present paper. And indeed, also for genuine
$2\times 2$ operator matrices, the blocks are often not assumed to
be closed as the closedness of at least some of the blocks would
be pointless anyway.

Henceforth, for $i=1,2$ let $E_i$ and $F_i$ be locally convex
Hausdorff spaces and for $i,j=1,2$ let $A_{ij}\in \mathfrak
L(E_{j},F_{i})$, where $\mathfrak L(E_{j},F_{i})$ denotes the set
of all densely defined operators from $E_j$ to $F_{i}$. The closed
densely defined and continuous everywhere defined operators in
these spaces are denoted by $\mathfrak C(E_{j},F_{i})$ and
$\mathfrak B(E_{j},F_{i})$, respectively. Finite direct sums will
be denoted with the symbol $\oplus$, meaning that in case of
Hilbert spaces, the direct sum will always be identified with the
orthogonal one. Define
  \[A=
  \begin{pmatrix}
    A_{11}&A_{12}\\A_{21}&A_{22}
  \end{pmatrix}
\]
from $E$ to $F$, where $E=E_{1}\oplus E_{2}$ and $F=F_1\oplus
F_2$, and the domain of $A$ is given by
   \begin{equation} \label{12.2}
   \mathcal D(A)=\left(\mathcal D(A_{11})\cap \mathcal
   D(A_{21})\right)\oplus \left(\mathcal D(A_{12})\cap \mathcal
   D(A_{22})\right).
   \end{equation}
The range of $A$ will be denoted by $\mathcal R(A)$.

Here we will only consider the case that $A$ is densely defined,
so we require that
   \begin{equation} \label{12.1}
  \text{$\mathcal D(A_{1j})\cap \mathcal D(A_{2j})$ is dense in
  $E_{j}$ for $j=1,2$.}
   \end{equation}

In particular, since all $A_{ij}$ are densely defined, their
adjoints $A_{ij}'$ must necessarily be closed operators from
$F_i'$ to $E_j'$, and we can define the operator
  \[A^\times =
  \begin{pmatrix}
    A_{11}'&A_{21}'\\A_{12}'&A_{22}'
  \end{pmatrix}
\]
from $F'$ to $E'$, where, according to \eqref{12.2}, $$\mathcal
D(A^\times )=\left(\mathcal D(A_{11}')\cap \mathcal D(A_{12}')
\right)\oplus \left(\mathcal D(A_{21}')\cap \mathcal
D(A_{22}')\right).$$

{Row and column operators of arbitrary size are defined as
follows.} A row operator $R\okr R_{R_1,\dots,R_n}$ is a linear
mapping defined in
  $\bigoplus_{j=1}^n E_j$, {taking values in $E_0$ and acting according to}
 \[R(\oplus _{j=1}^nf_j)=
\sum\nolimits_{j=1}^n R_jf_j,\quad
    \oplus _{j=1}^nf_j\in \mathcal
    D(R)\okr\bigoplus\nolimits_{j=1}^n \mathcal D(R_j),\]  where
    $E_0,E_1,\dots,E_n$ are locally convex Hausdorff spaces, and
    the $R_j$ are linear operators from $E_j$ to $E_0$. Clearly,
    $R$ is densely defined if and only if all $R_j$ are densely
    defined.

A column operator $C\okr C_{C_1,\dots,C_n}$ is a linear mapping
defined in $F_0${, taking values in $\bigoplus_{j=1}^n F_j$ and
acting as}
\[\mathcal D(C)\okr \bigcap\nolimits_{j=1}^n D(C_j),\quad Cf\okr \bigoplus\nolimits_{j=1}^n C_jf,\quad
    f\in \mathcal D(C),\]
 the
spaces $F_0,F_1,\dots,F_n$ are locally convex Hausdorff spaces,
and $C_j\in \mathfrak L(F_0,F_j)$ for $j=1,\dots,n$. To such a
column operator we associate the row operator
  \[C^\times \okr R_{C_1',\dots,C_n'}\]
  which may not be densely defined even if $C$ is densely defined.

   The important notice we can make at this stage more precise is that row
   operators behave differently than the column ones. In
   particular, supposing all the entries of $R$ are closed,
   \begin{equation*}
   \text{$R''$ may not be equal to $R$ though always
   $R'^\times=R$.}
   \end{equation*}
   Needless to say a similar behaviour concerns column operators.
   This outlines {once more} the flavour of our paper.

In Hilbert spaces, we will take the usual Hilbert space adjoints
of operators, i.\,e., with respect to the sesquilinear scalar
product rather than the adjoints with respect to bilinear forms
for dual pairs. Since the results of this paper are independent of
whether the duality is realized by bilinear forms or sesquilinear
forms, in Hilbert space we will always use Hilbert space adjoints.
With some abuse of notation, we shall write e.\,g.{} $C^\times $
in both cases, where its meaning will be clear from the context.
This is {applicable in particular to} Kre\u \i n space adjoints if each of the
component spaces is a Kre\u \i n space because in this case the
fundamental symmetry on the direct sum is the direct sum of the
fundamental symmetries on the components.

\section{Row operators}
In this section let $R\okr R_{R_1,\dots,R_n}$ be a row operator as
defined above.

\begin{prop}\label{row:1}\footnote{\;Propositions \ref{row:1},
\ref{col:1} and \ref{Cbounded} have been proved for linear
relations in Hilbert spaces in \cite{henk}, Proposition 2.1.} The
adjoint $R'$ of a densely defined row operator $R$ is the column
operator formed by the adjoints of the $R_j$, that is
$$R'=C_{R_1',\dots,R_n'}.$$
\end{prop}

\begin{proof}
To show $R'\subset C_{R_1',\dots,R_n'}$ let all $f_j\in \mathcal
D(R_j)$ and $g\in \mathcal D(R')$. Then
  \begin{equation*}
    \langle \oplus_j\! f_j, R'g\rangle
   =\big\langle\sum\nolimits_j R_jf_j,g\big\rangle
   =\sum\nolimits_j\langle R_jf_j,g\rangle.
  \end{equation*}
Putting $f_j=0$ for $j\ne k$, it follows that $g\in \mathcal
D(R_k')$ for each $k=1,\dots,n$ and
  $\langle R_kf_k,g\rangle =\langle f_k,R_k'g\rangle $.
Thus
  \[    \langle \oplus_j\! f_j, R'g\rangle
    =\langle \oplus_j\! f_j, \oplus_j\! R_j'g
    \rangle ,\] which proves the inclusion ``$\subset $''.

Conversely, let $g\in \mathcal D \left(C_{R_1',\dots,R_n'}\right)
=\bigcap_{j=1}^n\mathcal D(R_j')$. Then, for all $f_j\in \mathcal
D(R_j)$, $\langle f_j,R_j'g\rangle =\langle R_jf_j,g\rangle ,$ and
thus
  \begin{equation*}
    \langle \oplus_j f_j,\oplus_j R_j'g\rangle
   =\big\langle \sum\nolimits_j R_jf_j,g\big\rangle \\
   =\langle R\big(\oplus_j f_j\big),g\rangle,
  \end{equation*}
which shows $g\in \mathcal D(R')$.
\end{proof}
   \begin{rem} \label{tx1}
   The simplest example showing that a row operator with closed
entries may not be even closable is to consider $R_1$ and $R_2$
selfadjoint with $\mathcal D(R_1)\cap\mathcal
D{(R_2)}=\varnothing$. Then, according to Proposition \ref{row:1}
$\mathcal D{(R')}=\mathcal{D}{(R_1')}\cap\mathcal
D{(R_2')}=\varnothing$, hence $R$ is not closable. This is rather
an extreme example in a sense, a richer one concerning the same
question will be given below.
   \end{rem}

\begin{exa}\label{E:1}
Define $R_1$ and $R_2$ in $\ell^2$ as follows:
  \begin{equation*}
    R_1e_k\okr k^2e_k,\quad R_2e_k\okr ke_0+k^2e_k,\quad
    k=0,1,2,\dots,
  \end{equation*}
 where $(e_n)_{n=0}^\infty$ is the orthodox zero-one orthonormal
basis of $\ell^2$. To properly establish $R_1$ and $R_2$, we
define their domains by
  \begin{align*}
    \mathcal D(R_1)&=\left\{f=\sum _{k=0}^\infty \gamma
    _ke_k:
    \sum_{k=0}^\infty |\gamma _kk^2|^2<\infty
  \right\},\\
    \mathcal D(R_2)&=\left\{f=\sum _{k=0}^\infty \gamma _ke_k\in
    \mathcal D(R_1): \sum_{k=0}^\infty \gamma _kk\text{
    converges}\right\}.
  \end{align*}
Then, as usual, the operators $R_i$ are defined as
  \[R_if=\sum_{k=0}^\infty \gamma _kR_ie_k,\quad f=\sum_{k=0}^\infty \gamma _ke_k\in \mathcal D(R_i).\]
 Because $R_1$ is a diagonal operator on its maximal domain, it is
closed. In order to prove the closedness of $R_2$ let
  \[ f_n=\sum _{k=0}^\infty \gamma _{nk}e_k\in \mathcal D(R_2),\quad
f_n\to f=\sum _{k=0}^\infty \gamma _ke_k,\text{
and}\quad
  R_2f_n\to g=\sum_ {k=0}^\infty \delta _ke_k.\] If $P$ denotes
the orthogonal projection of the Hilbert space onto the closed
linear span of $\{e_k\}_{k=1}^\infty $, then $PR_2f_n=R_1f_n$, and
so $R_1f_n\to Pg$. Since $R_1$ is closed, it follows that $f\in
\mathcal D(R_1)$ and $Pg=R_1f$.

    Note that
  \begin{equation}\label{row2}
   \sum_{k=1}^\infty |\gamma _{nk}k^2-\gamma _kk^2|^2
   =\|R_1f_n-R_1f\|^2\to 0\quad\text{as}\quad n\to\infty,
  \end{equation}
  \[(I-P)R_2f_n=\sum_{k=1}^\infty \gamma _{nk}ke_0\]
and that $R_2f_n\to g$ implies
  \begin{equation}\label{row3}
    \sum _{k=1}^\infty \gamma _{nk}k\to\delta
_1\quad\text{as}\quad n\to\infty.
  \end{equation}
 Because for $m\ge 1$ and $n=0,1,\dots$,
  \begin{gather*}
   \big|\sum_{k=0}^m\gamma _kk-\delta _1\big|\le
   \big|\sum_{k=0}^m\gamma _{nk}k-\delta _1\big| +
   \sum_{k=0}^m|\gamma _{nk}k-\gamma _kk|\\\le
   \big|\sum_{k=0}^m\gamma _{nk}k-\delta _1\big| +
   \big(\sum_{k=1}^\infty |\gamma _{nk}k^2-\gamma _kk^2
   |^2\big)^{\frac12}\big(\sum_{k=1}^\infty \big|k
   \big|^{-2}\big)^{\frac12},
  \end{gather*}
from \eqref{row2}, and \eqref{row3} we get that for each
$\varepsilon >0$ there is $n\in \mathbb N$ such that
  \begin{equation}\label{row4}
  \big|\sum_{k=0}^m\gamma _{k}k-\delta _1\big|\le
  \big|\sum_{k=m+1}^\infty \gamma _{nk}k\big|+\varepsilon
  \end{equation}
for $m\ge 1$. From \eqref{row4} and the fact that $\sum
_{k=1}^\infty \gamma _{nk}k$ converges, we deduce that
  \[ \big|\sum_{k=1}^m\gamma _{k}k-\delta_1\big|\le2\varepsilon\]
for sufficiently large $m$, which proves that $\sum_{k=1}^\infty
\gamma _kk$ converges to $\delta _1$. This shows $f\in
\mathcal D(R_2)$ and $g=R_2f$.

 Finally, letting for $n\ge 1$ $f_n=n^{-1}e_n$, we have $f_n\in
\mathcal D(R_2)$, $f_n\to 0$ as $n\to \infty $, and
  \[ R_1(-f_n)+R_2f_n=nn^{-1}e_0=e_0.\]
So $R$ is not closable. However the linear span of
$(e_n)_{n=0}^\infty$ is included in $\mathcal D(R_1)\cap \mathcal
D(R_2)$ as well as in $\mathcal D(R_1')\cap \mathcal D(R_2')$ and
this is what makes this example more interesting than that argued
for in Remark \ref{tx1}.
\end{exa}

\begin{cor}\label{row:2}
For $R$ to be closable it is necessary but not sufficient that
$R_1, \dots,R_n$ are closable.
\end{cor}

\begin{proof} Use all these above and the fact, which is implicit in
Proposition \ref{row:1}, that $\mathcal
D(R')=\bigcap_{j=1}^n\mathcal D(R_j')$.
\end{proof}

Note that in general $R$ is not closed even if all its entries are
closed. However, we trivially have

\begin{prop}\label{Rbounded1}
Assume that at most one of the entries of $R$ is not bounded. Then
$R$ is closed (closable) if and only if all $R_j$ are closed
(closable). In this case $R'$ is densely defined.
\end{prop}

 In the following we have a particular result when the closure of
$R$ can be determined explicitly.

\begin{prop}\label{row:3}
Let $n=2$, assume that $R_1$ is injective, that $\mathcal
R(R_2)\subset \mathcal R(R_1)$ and that $R_1^{-1}R_2$ is an
operator with a bounded extension $K\in \mathfrak B(E_2,E_1)$.
Then
  \[\overline{R}=(\,\overline{R_1},0)
  \begin{pmatrix}
    I&K\\0&I
  \end{pmatrix},\]
and $\mathcal D(\,\overline{R}\,)= \{f\oplus g\in E_1\oplus
E_2:f+Kg\in \mathcal D(\,\overline{R_1}\,)\}$.
\end{prop}

\begin{proof}
The operator $\hat R=(\,\overline{R_1},0)
  \begin{pmatrix}
    I&K\\0&I
  \end{pmatrix}$ is closed since the $2\times 2$ matrix operator on the
right is invertible. Also, for $f\in \mathcal D(R_1)$, $g\in
\mathcal D(R_2)$ we have $Kg=R_1^{-1}R_2 g\in \mathcal D(R_1)$ and
thus $f+Kg\in \mathcal D(R_1)$, and
  \begin{equation*}
    \hat R (f\oplus g) =(\,\overline{R_1},0) (f+R^{-1}R_2g)\oplus
    g =R_1f+R_2g=R (f\oplus g).
  \end{equation*}
This shows that $R\subset \hat R$ and thus $\overline{R}\subset
\hat R$ as $\hat R$ is closed. We calculate
  \begin{equation}\label{row5}
\hat R'=
  \begin{pmatrix}
    I&0\\K'&I
  \end{pmatrix}
  \begin{pmatrix}
    R_1'\\0
  \end{pmatrix}=
  \begin{pmatrix}
    R_1'\\K'R_1'
  \end{pmatrix}.
  \end{equation}
From $R_2\subset R_1K$ it follows that
  $K'R_1'\subset (R_1K)'\subset R_2'$.
In particular, $\mathcal D(R_1')\subset \mathcal D(R_2')$ and
hence $K'R_1'=R_2'$ on $\mathcal D(R_1')$. Thus $\hat R'=R'$ by
\eqref{row5} and Proposition \ref{row:1}, which proves $\hat
R=\overline{R}$ since $\hat R$ is closed. The representation of
$\mathcal D(\,\overline{R}\,)$ is now obvious from the
representation of $\overline{R}$.
\end{proof}

\section{column operators}

In this section, let $C\okr C_{C_1,\dots,C_n}$ be a densely
defined column operator.

\begin{prop}\label{col:1}
$C^\times \subset C'$.
\end{prop}

\begin{proof}
Since $C$ is densely defined, $C'$ is an operator. Let $g_j\in
\mathcal D(C_j')$,  $f\in \mathcal
D(C)$. Then
  \begin{align*}
    \langle f,C^\times (\oplus_j g_j)
\rangle =\langle f,\sum\nolimits_j C_j'g_j\rangle    =\sum\nolimits_j\langle f,C_j'g_j\rangle
   =\sum\nolimits_j\langle C_jf,g_j\rangle
    =\big\langle Cf,\oplus_j g_j\big\rangle .
  \end{align*}
This proves $\oplus_j g_j\in \mathcal D(C')$ and
$C'(\oplus_j g_j)
=C^\times(\oplus_j g_j)$.
\end{proof}

   Referring back to the preceding section, let us remind that
$R^\times= C_{R_1',\dots,R_n'}$. Then the conclusion of
Proposition \ref{row:1} can be restated as
   \begin{equation*}
   R'=R^\times.
   \end{equation*}
    Now double applications of this and the Proposition
\ref{col:1} leads to
   \begin{cor} \label{ty1}
   For a row operator $R$ and a column operator $C$ we have
   \begin{equation*}
   R''=R^{\times\prime}\supset
R^{\times\times}=R^{\prime\times},\quad C''\subset
C^{\times\prime}=C^{\times\times}.
   \end{equation*}
   \end{cor}

   \begin{prop}\label{col:4}
\liczp 1 $C$ is closed if $C_1,\dots,C_n$ are closed. \liczp 2 $C$
is closable if $C_1,\dots,C_n$ are closable. \liczp 3
$\overline{C^\times }=C'$ if and only if $\overline{C}=
C_{\,\overline{ C_1},\dots,\overline{C_n}\,}$.
\end{prop}

\begin{proof}
\liczp 1 is straightforward, and \liczp 2 is an immediate
consequence of \liczp 1. From Corollary \ref{ty1} we get
   \begin{equation*}
   \overline{C^\times}=C^\prime \Longleftrightarrow
C^{\times\prime}=C^{\prime\prime} \Longleftrightarrow
C^{\times\times} =C^{\prime\prime}.
   \end{equation*}
   Then the conclusion follows by observing that
$C^{\prime\prime}=\overline{C}$ and
   \begin{equation*}
   C^{\times\times}=R_{C_1^\prime,\dots,C_n^\prime}^\times=
C_{\overline{C_1},\dots,\overline{C_n}}.
   \end{equation*}
   \end{proof}

\begin{rem}\label{col:3} The sufficient condition for $C$ to be closable, which is in
Proposition \ref{col:4}, \liczp 2, turns out to be not necessary.
For this let $H_1$ and $H_2$ be Hilbert spaces and $C\in\mathfrak
C(H_1, H_2) \setminus \mathfrak B(H_1,H_2)$, that is, $\mathcal
D(C^*)\ne H_2$. Let $x\in H_2\setminus \mathcal D(C^*)$, $P$ be
the orthogonal projection onto the span of $x$, $C_1=PC$,
$C_2=(I-P)C$. Then $C=
  \begin{pmatrix}
    C_1\\C_2
  \end{pmatrix}$
with $C_1$ not closable. Here, with some abuse of notation, the
range spaces of $C_1$ and $C_2$ are $\mathcal R(P)$ and $\mathcal
R(I-P)$, respectively.

Indeed, assume that $x\in \mathcal D(C_1^*)$. Then we have for
$f\in \mathcal D(C)$ that
  \[\langle Cf,x\rangle =\langle Cf,Px\rangle =\langle PCf,x\rangle
    =\langle f,C_1^*x\rangle ,\] and the contradiction $x\in
    \mathcal D(C^*)$ would follow. Since $\mathcal D(C_1^*)\subset
    \mathcal R(P)$ and $\mathcal R(P)$ is one-dimensional,
    $\mathcal D(C_1^*)=\{0\}$ follows. Thus $C_1$ is not closable.

For other examples with nonclosable $C_1$ but closable $C_2$ and
$C$ we refer to \cite[Sections 2 and 3]{DGS} and \cite[Section
1]{N}.
\end{rem}

\begin{prop}\label{Cbounded}
Assume that at most one of the entries of $C$ does not satisfy
$C_j\in\mathfrak B(F_0,F_j)$. Then $C^\times =C'$.
\end{prop}

\begin{proof}
Obviously, $\overline{C}=
C_{\,\overline{C_1},\dots,\overline{C_2}}$, and hence
$\overline{C^\times }=C'$ by Proposition \ref{col:4} \liczp 3.
 Clearly, $C^\times $ satisfies the assumptions of Proposition
\ref{Rbounded1}, and thus $C^\times$ is closed.
\end{proof}

 In the following result, $C^\times $ and $C'$ are given
explicitly; in particular, the structure of their domains is
transparent.

\begin{prop}\label{col:5}
Let $C_1$, $C_2$ be such that $\mathcal D(C_1)=\mathcal D(C_2)$
and assume that $C_1$ is injective with $C_2C_1^{-1}$ having an
extension $K\in \mathfrak B(F_1,F_2)$. Then
\[ C^\times =(C_1',C_1'K') \quad \text{and}\quad
  C'=(C_1',0)
  \begin{pmatrix}
    I&K'\\0&I
  \end{pmatrix}.\]
In particular, $\mathcal D(C^\times ) =\{f\oplus g\in F_1'\oplus
F_2':f\in \mathcal D(C_1'),\, K'g\in \mathcal D(C_1')\}$ and
$\mathcal D(C') =\{f\oplus g\in F_1'\oplus F_2':f+K'g\in \mathcal
D(C_1')\}$.
\end{prop}

\begin{proof}
Since $C_2=C_2C_1^{-1}C_1=KC_1$, $C_2'=C_1'K'$, and the statements
about $C^\times $ follow.

For $C$ we have the representation
  \[C=
  \begin{pmatrix}
    C_1\\KC_1
  \end{pmatrix}=
  \begin{pmatrix}
    I&0\\K&I
  \end{pmatrix}
  \begin{pmatrix}
    C_1\\0
  \end{pmatrix},\]
which gives
   \[C'=(C_1',0)
  \begin{pmatrix}
    I&K'\\0&I
  \end{pmatrix}.\]
The statement about the domain of $C'$ immediately follows from
this.
\end{proof}

   \section{More about column operators in Hilbert spaces}
Here we assume that $F_0$, $F_1$, $F_2$ are Hilbert spaces and
that $C=C_{C_1,C_2} $ is a densely defined closable column
operator.
Let $D_0$ be a subspace of $F_0$ satisfying
  \begin{equation}\label{colH2}
    D_0\subset \mathcal D(\overline{C}C^\times ).
  \end{equation}
   Notice that, if $\overline{C^\times}=C^*$ the, by part \liczp 3
of Proposition {\rm\ref{col:4}},
  \[\overline{C}C^\times =
  \begin{pmatrix}
    \overline{C_1}C_1^*&\overline{C_1}C_2^*\\
   \overline{C_2}C_1^*&\overline{C_2}C_2^*
  \end{pmatrix}\]

We want to investigate the following
questions\,\footnote{\;$\mathcal D\subset\mathcal D(A)$ is a core
for a closable operator $A$ if $A\subset \overline{A|_{\mathcal
D}}$.}:
\begin{itemize}
\item when is $D_0$ a core for $C^*$, \item when is $D_0$ a core
for $C^\times $?
\end{itemize}

\begin{prop}\label{colH:1}
Let $D_1$ be a subspace of $ F_0$ such that $D_0\subset D_1\subset
\mathcal D(C^*)$. Then $D_0$ is a core for $C^*|_{D_1}$ if and
only if
   \begin{equation} \label{XXX}
   ((I+\overline{C}C^\times )(D_0))^\perp\cap D_1=\{0\}.
   \end{equation}
   \end{prop}

\begin{proof}
Due to \eqref{colH2}, for $f\in D_1$ to belong to the left hand
side of \eqref{XXX} means precisely that \begin{equation*}
   \text{$0= \langle f,(I+\overline{C}C^*)g\rangle$ for all $g\in
D_0$.}
  \end{equation*}
   Because $D_1\subset \mathcal D(C^*)$, this is what is required
for $D_1$ to be a core for $C^*$.
\end{proof}

\begin{prop}\label{colH:2}
Let $C_1$, $C_2$ be such that $\mathcal D(C_1)=\mathcal D(C_2)$
and assume that $C_1$ is injective with $C_2C_1^{-1}$ having an
extension $K\in \mathfrak B(F_1,F_2)$. Let $D_0$ be a subspace of $\mathcal
D(C^*)$. Also assume that there is at least one $v_0\in
F_2\setminus \{0\}$ such that $(-K^*v_0,v_0)\in D_0$. Then $D_0$
is a core for $C^*$ if and only if $D_0$ is dense in $F_1\oplus
F_2$ and $(I,K^*)(D_0)\cap\mathcal D(C_1^*)$ is a core for
$C_1^*$.
\end{prop}

\begin{proof}
Since we assume that $C$ is closable, $C^*$ is densely defined,
and the condition that $ D_0$ is dense in $F_1\oplus F_2$ is
necessary for a core. Thus we may assume this property for the
remainder of the proof. Using the graph norm of $C^*$ as in the
proof of Proposition \ref{colH:1}, we have to investigate when,
for $f\oplus g\in \mathcal D(C^*)$,
  \begin{equation}\label{colH4}
    \langle u\oplus v,f\oplus g\rangle +\langle C^*(u\oplus v),C^*
    (f\oplus g)\rangle =0
  \end{equation}
for all $u\oplus v\in D_0$ implies $f\oplus g=0$. By Proposition
\ref{col:5},
   \[C^*(u\oplus v)=C_1^*(u+K^*v) \quad\text{and}\quad
    C^*(f\oplus g)=C_1^*(f+K^*g).\] Introducing $w=u+K^*v$ and
    $h=f+K^*g$, \eqref{colH4} can be written as
  \begin{equation}\label{colH5}
    \langle u,f\rangle +\langle v,g\rangle +\langle
    C_1^*w,C_1^*h\rangle =0,
  \end{equation}
which is the same as
  \begin{equation*}
    \langle w-K^*v,h-K^*g\rangle +\langle v,g\rangle +\langle
    C_1^*w,C_1^*h\rangle =0.
  \end{equation*}
Note that
   \[  \langle w-K^*v,h-K^*g\rangle +\langle v,g\rangle -\langle w,h\rangle
    =\langle (I+KK^*)v-Kw,g\rangle -\langle K^*v,h\rangle .\]
    First we want to find a condition for $h=0$; thus we may
    modify $f$ and $g$ as long as $f+K^*g=h$. If $(I+KK^*)v\ne Kw$
    choose $g\in F_2$ such that
  \[\langle (I+KK^*)v-Kw,g\rangle =\langle K^*v,h\rangle .\]
If $(I+KK^*)v= Kw$, we replace $u\oplus v$ with $(u-K^*v_0)\oplus
(v+v_0)$. Since $(u-K^*v_0)+K^*(v+v_0)=u+K^*v=w$, this does not
change $w$, but now
   \[(I+KK^*)(v+v_0)=Kw+(I+KK^*)v_0\ne Kw\]
since $v_0\ne0$ and $(I+KK^*)$ is injective. Thus also in this
case we can find $g\in F_2$ such that
  \[\langle (I+KK^*)(v+v_0)-Kw,g\rangle =\langle K^*(v+v_0),h\rangle .\]
With this $g$ we put $f=h-K^*g$ in either case and deduce
   \[\langle w,h\rangle +\langle C_1^*w,C_1^*h\rangle =0\]
for all $w\in (I,K^*)(D_0)$. Hence $h =0$ for all $h\in \mathcal
D(C_1^*)$ if and only if $(I,K^*)(D_0)\cap\mathcal D(C_1^*)$ is a
core for $C_1^*$. Returning to \eqref{colH5} it follows (now for
our original $f$ and $g$) that $h=0$ implies
  \[\langle u,f\rangle +\langle v,g\rangle =0\]
for all $u\oplus v\in D_0$. Since $D_0$ is dense in $F_1 \oplus
F_2$, $f=0$ and $g=0$ follow.
\end{proof}

\begin{cor}\label{colH:3}
Let $C_1$, $C_2$ be such that $\mathcal D(C_1)=\mathcal D(C_2)$
and assume that $C_1$ is injective with $C_2C_1^{-1}$ having an
extension $K\in \mathfrak B(F_1,F_2)$. Then $\overline{C^\times }=C^*$.
\end{cor}

\begin{proof}
For every $v_0\in \mathcal D(C_1^*K^*)$ we have $-K^*v_0\in
\mathcal D(C_1^*)$, so that $(-K^*v_0,v_0)\in \mathcal D(C^\times
)$ holds for every (and thus some) $v_0\in F_2\setminus \{0\}$
since $\mathcal D(C^\times )=\mathcal D(C_1^*)\oplus \mathcal
D(C_1^*K^*)$ by Proposition \ref{col:5}. Obviously,
$(I,K^*)(\mathcal D(C^\times )) \supset \mathcal D(C_1^*)$, and an
application of Proposition \ref{colH:2} completes the proof.
\end{proof}

   \begin{exa} \label{t15.9.1}
   Suppose $C$ is a column operator with $C_2=T-C_1$, both $C_1$
   and $T$ are in $\llf(F,F_0)$, where $F\okr F_1=F_2$. Suppose
   moreover
   \begin{equation*}
   \text{$\dz{C_1}\subset\dz{T}$ and
   $\dz{C_1\gw}\varsubsetneq\dz{T\gw}$.}
   \end{equation*}
Then
\begin{align*}
   \dz{C^\times}\varsubsetneq \dz{C\gw}.
   \end{align*}
Indeed, for $h\in\mathcal D(C_1), \;f,g\in F$
   \begin{equation*}
   \is{Ch}{f\oplus g}=\is{C_1h\oplus(T-C_1)h}{f\oplus
   g}=\is{C_1h}{f-g}+\is{Th}g.
   \end{equation*}
   This gives us immediately that $f\oplus f$ is in $\dz{C\gw}$ if
   and only if $f$ is in $\dz{T\gw}$; then
   \begin{equation*}
   C\gw (f\oplus f)=T\gw f,\quad f \in\dz{T\gw}
   \end{equation*}
    and consequently
   \begin{align*}
   \dz{C^\times}\subset\dz{C^\times}+\zb{f\oplus
   f}{f\in\dz{T\gw}}\subset\dz{C\gw}.
   \end{align*}
    Hence, if $f\in\dz{T\gw}\setminus\dz{C_1\gw}$, then $f\oplus
    f\in \dz{C\gw}$ but $f\oplus
    f\notin\dz{C_1\gw}\bigoplus \dz{C_2\gw}=\dz{C^\times}$.

   \end{exa}

\section{The operator matrix $A$}

Let $A$ be a $2\times 2$ matrix of the form
   \begin{equation*}
  \begin{pmatrix}
    A_{11}&A_{12}\\A_{21}&A_{22}
  \end{pmatrix}
   \end{equation*}
 with the denseness condition \eqref{12.1} being satisfied.
Writing
  \begin{equation}\label{A0}
    C_{A,i}=(A_{i1},A_{i2}),\ i=1,2,
  \end{equation}
we know from Proposition \ref{row:1} that
  \begin{equation*}
 C_{A,i}'=
  \begin{pmatrix}
    A_{i1}'\\A_{i2}'
  \end{pmatrix},
  \end{equation*}
and thus
  \begin{equation}\label{A1}
    A=
  \begin{pmatrix}
    C_{A,1}\\C_{A,2}
  \end{pmatrix},\quad A^\times =(C_{A,1}',C_{A,2}').
  \end{equation}

\begin{thm}\label{main}
\liczp 1 $A^\times \subset A'$.
\liczp 2 $A^\times$ is closable.
\liczp 3 If $A^\times$ is densely defined, then $A$ is closable.
\end{thm}
\begin{proof}
\liczp 1 immediately follows from \eqref{A1} and Proposition \ref{col:1}.  Since $A$ is densely defined, $A'$ is a closed operator, and
hence by \liczp 1,
$A^\times $ is closable, which makes \liczp 2.
If $A^\times $ is densely defined, so is $A'$ by \liczp 1. Consequently, $A$ is
closable and \liczp 3 follows.
\end{proof}

\begin{prop}\label{ret:1}
If $\mathcal D(C_{A,1}')\oplus \mathcal D(C_{A,2}')$ is dense in
$F'$, then $A$ is closable.
\end{prop}

\begin{proof}
If $\mathcal D(C_{A,1}')\oplus \mathcal D(C_{A,2}')$ is dense in
$F'$, then $A'$ is densely defined by part \liczp 1 of Theorem
\ref{main}, which means that $A$ is closable.
\end{proof}

More can be said if at most one of the operators $A_{ij}$ is not
bounded.

\begin{prop}\label{Abdd}
Assume that $A_{ij}\in \mathfrak B(E_j,F_i)$ with the exception of
at most one pair $(i,j)$ and that this exceptional $A_{ij}$ is
closable. Then $A$ is closable and $A'=A^\times $.
\end{prop}

\begin{proof}
By Proposition \ref{Rbounded1}, both $C_{A,1}$ and $C_{A,2}$ are
closable, and hence $A$ is closable by Proposition \ref{col:4} \liczp2.
 At most one of the operators $C_{A,i}$ in \eqref{A0} does not
satisfy $C_{A,i}\in \mathfrak B(E,F_i)$. Then Proposition
\ref{Cbounded} and \eqref{A1} lead to $A'=A^\times $.
\end{proof}

Theorem \ref{main} and Proposition \ref{Abdd} immediately raise
the question if the following cases can occur:

\begin{enumerate}
\renewcommand{\labelenumi}{\Roman{enumi}.}
\item A is not closable, \item $A$ is closable and $A^\times $ is
not densely defined, \item $A^\times $ is densely defined but
$A'\ne \overline{A^\times }$, \item $A'=\overline{A^\times }$ but
$A'\ne A^\times $.
\end{enumerate}

Below we will show indeed that all these cases can occur, even
under the additional requirement that all $A_{ij}$ are closed
operators in Hilbert spaces.
\begin{exa}\label{I} Here we give an example for I. Let $E_1=E_2=F_1=F_2
=\ell^2$ and, with $R_1$, $R_2$ from Example \ref{E:1}, put
$A_{11}=R_1$, $A_{12}=R_2$, $A_{21}=A_{22}=0$. Then all $A_{ij}$
are closed, and since $(R_1,R_2)$ is not closable, also $A$ is not
closable by Proposition \ref{col:4} \liczp 2.
\end{exa}

\begin{exa}\label{II} Here we give an example for II. Let $E_1=E_2=F_1=F_2
=\ell^2$ and, with $R_1$, $R_2$ from Example \ref{E:1}, put
$A_{11}=R_1$, $A_{12}=R_2$, $A_{21}=0$, $A_{22}=R_2$. Then all
$A_{ij}$ are closed, $A$ is closable, {while} $A^\times $ is not
densely defined.
\end{exa}

\begin{proof}
To show that $A$ is closable, consider any sequences
  \[ f_n\in \mathcal D(A_{11})\cap \mathcal D(A_{21}),\
     g_n\in \mathcal D(A_{12})\cap \mathcal D(A_{22})\] satisfying
     $f_n\to0$, $g_n\to 0$, $A_{11}f_n+A_{12}g_n\to h_1$ and
     $A_{21}f_n+A_{22}g_n\to h_2$ for some $h_1$, $h_2$ in
     $\ell^2$. Then
  \[h_2=\lim_{n\to\infty }(A_{21}f_n+A_{22}g_n)=\lim_{n\to\infty }R_2g_n,\]
$g_n\to0$, and the closedness of $R_2$ imply that $h_2=0$.
Consequently,
  \[h_1=\lim_{n\to\infty }(A_{11}f_n+A_{12}g_n)
    =\lim_{n\to\infty }(R_1f_n+R_2g_n)=\lim_{n\to\infty }R_1f_n,\]
    so that $f_n\to 0$ and the closedness of $R_1$ imply $h_1=0$.
    This completes the proof of the closability of $A$. By Example
    \ref{E:1}, $(A_{11},A_{12})$ is not closable, whence, in view
    of Proposition \ref{row:1}, $
  \begin{pmatrix}
    A_{11}'\\A_{12}'
  \end{pmatrix}$
is not densely defined. Hence also $A^\times $ is not densely
defined.
\end{proof}

\begin{exa}\label{III}
Here we give an example for III. Let $E_1=E_2=F_1=F_2=L_2(0,1)$
and let the operators $A_{ij}$, $i,j=1,2$ be defined by
  \begin{align*}
    \mathcal D(A_{11})&=\{f\in W_2^1(0,1):f(0)=0\}, \ A_{11}f=f',\\
     \mathcal D(A_{12})&=L_2(0,1),\ A_{ 12}=0,\\
     \mathcal D(A_{21})&=\{f\in W_2^1(0,1):f(1)=0\},\ A_{21}f=f',\\
     A_{22}&=-A_{21},
  \end{align*}
where $W_2^1(0,1)$ denotes the usual Sobolev space of order $1$.
Then $A$ is densely defined, all $A_{ij}$, $i,j=1,2$, are closed,
$A^\times $ is densely defined, and $\overline{A^\times}\ne A'$.
\end{exa}

\begin{proof}
The denseness of the domain of $A$ as well as the closedness of
the $A_{ij}$ is well-known and obvious. Since $A_{11}^*=-A_{21}$,
it is also clear that $A^\times $ is densely defined. Putting
  $C_1=(A_{11},A_{12}),\ C_2=(A_{21},A_{22}),$
we have
    \[A=
  \begin{pmatrix}
    C_1\\C_2
  \end{pmatrix}.\]
We also let $C_{12}=C_1|_{\mathcal D(C_1)\cap \mathcal D(C_2)}$
and $C_{21}=C_2|_{\mathcal D(C_1)\cap \mathcal D(C_2)}$. Clearly,
   \[A=
  \begin{pmatrix}
    C_{12}\\C_{21}
  \end{pmatrix}\subset
  \left(\,\begin{matrix} \overline{C_{12}}\\[2pt]\overline{C_{21}}
    \end{matrix}\,\right)\subset
  \left(\,\begin{matrix} \overline{C_1}\\[2pt]\overline{C_2}
  \end{matrix}\,\right)
,\] and in view of Proposition \ref{col:4} \liczp 1 it follows
that
  \begin{equation}\label{III.1}
    \overline{A} \subset \left(\,\begin{matrix}
    \overline{C_{12}}\\[2pt]\overline{C_{21}}
    \end{matrix}\,\right)\subset
  \left(\,\begin{matrix} \overline{C_1}\\[2pt]\overline{C_2}
  \end{matrix}\,\right).
  \end{equation}
Hence, by Proposition \ref{col:4} \liczp 3, the proof will be
complete if we show that the second inclusion in \eqref{III.1} is
strict, i.\,e.,
  \begin{equation}\label{III.2}
    \mathcal D(\,\overline{C_{12}}\,)\cap\mathcal
    D(\,\overline{C_{21}}\,) \ne \mathcal
    D(\,\overline{C_{1}}\,)\cap\mathcal D(\,\overline{C_{2}}\,).
  \end{equation}
Since $C_2(f\oplus f)=0$ for $f\in \mathcal D(A_{21})$ we have
$\{f\oplus f:f\in L_2(0,1)\}\subset \mathcal
D(\,\overline{C_2}\,)$, which immediately leads to
  \begin{equation}\label{III.3}
    \{f\oplus f:f\in \mathcal D(A_{11})\}\subset \mathcal
    D(\,\overline{C_1}\,) \cap \mathcal D(\,\overline{C_2}\,).
  \end{equation}
Since $A_{11}$ and $A_{11}|_{\mathcal D(A_{11})\cap \mathcal
D(A_{21})}$ are closed operators, so are $C_1$ and $C_{12}$. Hence
$f\oplus g\in \mathcal D(\,\overline{C_{12}}\,)\cap \mathcal
D(\,\overline{C_{21}}\,)$ implies $f\in \mathcal D(A_{11})\cap
\mathcal D(A_{21})$, and thus \eqref{III.2} is proved in view of
\eqref{III.3} and $\mathcal D(A_{11})\not\subset \mathcal
D(A_{21})$.
\end{proof}

\begin{exa}\label{IV}
Here we give an example for IV. Let $E_1=E_2=F_1=F_2$ be an
infinite dimensional Hilbert space and let $A_{11}$ be a closed
densely defined unbounded operator in this Hilbert space. Put
$A_{12}=A_{22}=0$ and $A_{21}=-A_{11}$. Then $A$ is a densely
defined closed operator, $\overline{A^\times }=A^*$, and $A^\times
\ne A^*$.

Indeed, the operator $A$ is clearly densely defined and closed,
and
   \[  A^\times=
 \begin{pmatrix}
    A_{11}^* &-A_{11}^*\\0&0
  \end{pmatrix}.\]
  Considering $A^\times$ as a column operator of row operators
(one may look also at Theorem 10 in \cite{S}) and applying
Propositions \ref{Cbounded} and \ref{row:1} we come to
   \begin{equation*}
   A^{\times*}=\begin{pmatrix}
    A_{11}^{**} &0\\-A_{11}^{**}&0
  \end{pmatrix} =A
   \end{equation*}
    since $A_{11}$ is closed. Taking adjoints gives
$\overline{A^\times}=A^*$.

   Using Example \ref{t15.9.1} with $T=0$ we get $\mathcal
D(A^\times)\neq\mathcal D(A^*)$ which establishes $A^\times \ne
A^*$.

\end{exa}

\end{document}